\documentclass{segabs}

\usepackage{amsmath}
\usepackage{amsfonts}
\usepackage{bm} 

\newcommand{\bG}{\bm{\Gamma}}
\newcommand{\bhG}{\widehat{\bm{\Gamma}}}

\newtheorem{alg}{Algorithm}

\begin{document}

\title{S-fraction multiscale finite-volume method for spectrally accurate wave propagation}

\renewcommand{\thefootnote}{\fnsymbol{footnote}} 

\author{Vladimir Druskin, Alexander V. Mamonov\footnotemark[1] and Mikhail Zaslavsky, Schlumberger}

\footer{SFMSFV}
\lefthead{Druskin, Mamonov \& Zaslavsky}
\righthead{S-fraction multiscale finite-volume method}

\maketitle

\begin{abstract}
We develop a method for numerical time-domain wave propagation based on the 
model order reduction approach. The method is built with high-performance 
computing (HPC) implementation in mind that implies a high level of 
parallelism and greatly reduced communication  requirements compared to the 
traditional high-order finite-difference time-domain (FDTD) methods. The 
approach is inherently multiscale, with a reference fine grid model being 
split into subdomains. For each subdomain the coarse scale reduced order 
models (ROMs) are precomputed off-line in a parallel manner. The ROMs
approximate the Neumann-to-Dirichlet (NtD) maps with high (spectral) accuracy 
and are used to couple the adjacent subdomains on the shared boundaries. 
The on-line part of the method is an explicit time stepping with the coupled ROMs. 
To lower the on-line computation cost the reduced order spatial operator is
sparsified by transforming to a matrix Stieltjes continued fraction (S-fraction)
form. The on-line communication costs are also reduced due to the ROM NtD map 
approximation properties. Another source of performance improvement is the time 
step length. Properly chosen ROMs substantially improve the 
Courant-Friedrichs-Lewy (CFL) condition. This allows the CFL time step to 
approach the Nyquist limit, which is typically unattainable with traditional 
schemes that have the CFL time step much smaller than the Nyquist sampling rate.
\end{abstract}

\section{Introduction}

Seismic inversion requires simulating of acoustic or elastic wave propagation
on a very large scale. The computational cost of the forward problem is 
typically the dominant part of the overall inversion cost. Thus, fast and 
accurate wave propagation solvers capable of handling large models are of 
utmost importance. This is particularly true for the iterative inversion 
techniques such as the full waveform inversion (FWI), which require multiple 
forward solves at each iteration to compute the minimization functional and its 
derivative information. 

For maximum performance the method must be easily parallelizable. 
We achieve this by splitting the computational domain $\Omega$ into
the subdomains $\Omega^\alpha$. This is illustrated in Figure 
\ref{fig:scheme-bw-crop} for a 2D case. The derivations and the numerical 
examples below are for the full 3D case. 

On each subdomain $\Omega^\alpha$ we discretize the spatial operator on a 
fine reference grid. The resulting semi-discrete second order wave equation is
\begin{equation}
u_{tt}^\alpha = \tensor{A}^\alpha u^\alpha,
\label{eqn:wave}
\end{equation}
where $u^\alpha$ is the wave field at the fine grid nodes of the subdomain
$\Omega^\alpha$. Our framework provides a unified treatment for both the 
acoustic and elastic cases, so that we do not specify here the exact form of the 
fine grid discrete spatial operator 
$\tensor{A}^\alpha \in \mathbb{R}^{N \times N}$. The required boundary 
conditions are embedded into $\tensor{A}^\alpha$ as needed.

Our method consists of two stages. The first stage is performed off-line
before the time stepping. During this stage, the reduced order models are
computed for each $\tensor{A}^\alpha$. Because there is no interaction between
the subdomains at this point, the ROMs are computed in parallel. 

\plot{scheme-bw-crop}{width=0.55\columnwidth}
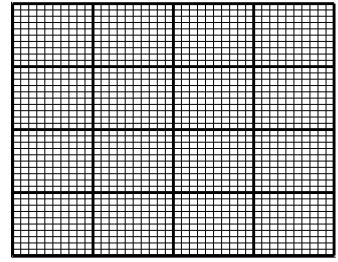
{Computational domain $\Omega$ with the reference fine grid (thin lines) 
split into $4 \times 4$ subdomains $\Omega^\alpha$ with the boundaries 
$\mathcal{B}^{\alpha\beta}$ (thick lines).}

The second stage is the time stepping. At this stage, the adjacent subdomains
exchange the information with each other at every time step. To maximize the 
overall performance, the ROMs constructed earlier should reduce the amount
of communication and make the time-step increase possible. 

The method presented here is an extension of the techniques of 
\cite[]{druskin2000gaussian, asvadurov2000application}, 
where the so-called optimal (spectrally matched) grids were used to construct
the ROMs on the subdomains. The use of optimal grids relies on the medium 
being uniform on each subdomain. The method presented here avoids this 
limitation and allows for arbitrary sharp discontinuities within the subdomains.

\section{Stage 1: reduced order models}

The two adjacent subdomains $\Omega^\alpha$ and $\Omega^\beta$ 
communicate only through the shared boundary $\mathcal{B}^{\alpha \beta}$. 
For a second order PDE all the exchanged information can be captured in a 
Neumann-to-Dirichlet (NtD) map. Thus, the ROM must approximate well the NtD 
map while reducing the number of degrees of freedom shared by the subdomains. 
To achieve this we choose first a small number $m$ of basis functions,
columns of $\tensor{F}^{\alpha\beta} \in \mathbb{R}^{N \times m} $, that 
are localized on $\mathcal{B}^{\alpha\beta}$ and are zero elsewhere in 
$\Omega^\alpha$. 

Let us denote by $\mathcal{N}(\alpha)$ the indices of the subdomains
adjacent to $\Omega^\alpha$. Then, we can combine all six sets 
(a 3D box has 6 faces) of basis functions 
$\tensor{F}^{\alpha\beta}, \; \beta \in \mathcal{N}(\alpha)$ into one 
matrix $\tensor{F}^\alpha \in \mathbb{R}^{N \times 6m}$. Transforming 
(\ref{eqn:wave}) to the frequency domain
\begin{equation}
\tensor{A}^\alpha u^\alpha + \omega^2 u^\alpha = 0,
\label{eqn:Helmholtz}
\end{equation}
we can write the frequency-dependent NtD map projected on the basis 
functions $\tensor{F}^\alpha$ as 
\begin{equation}
\tensor{M}^{\alpha}(\omega) = \left[ \tensor{F}^\alpha \right]^* 
\left( \tensor{A}^\alpha + \omega^2 \tensor{I} \right)^{-1} \tensor{F}^\alpha,
\label{eqn:NtD}
\end{equation}
which has the exact form of a transfer function of a multi-input/multi-output
(MIMO) dynamical system with both inputs and outputs given by 
$\tensor{F}^\alpha$. Hereafter we omit the subdomain index $\alpha$ unless
there are more than two subdomains under consideration at once.

Once the NtD map is expressed in the form (\ref{eqn:NtD}), we can apply the
well-developed theory of model order reduction to come up with a ROM
\begin{equation}
\widetilde{\tensor{M}} (\omega) = \widetilde{\tensor{F}}^* 
\left( \widetilde{\tensor{A}} + \omega^2 \tensor{I} \right)^{-1} \widetilde{\tensor{F}},
\label{eqn:NtDROM}
\end{equation}
where $\widetilde{\tensor{A}} \in \mathbb{R}^{6mn \times 6mn}$ and 
$\widetilde{\tensor{F}} \in \mathbb{R}^{6mn \times 6m}$ with $6mn \ll N$. 
The ROM has a block structure with $n$ being the number of blocks, as described
below.

To obtain high (spectral) accuracy of the resulting numerical scheme, we require
$\widetilde{\tensor{M}} (\omega)$ to be a good approximation of the NtD map 
$\tensor{M}(\omega)$ as a function of $\omega$. The existing
literature contains many approaches to this problem. A large family of 
approaches uses projection onto some subspace of $\mathbb{R}^N$ to obtain
(\ref{eqn:NtDROM}). If the columns of some 
$\tensor{V} \in \mathbb{R}^{N \times 6mn}$ form an orthonormal basis
for the desired projection subspace, the ROM is defined by
\begin{equation}
\widetilde{\tensor{A}} = \tensor{V}^* \tensor{A} \tensor{V}, \quad
\widetilde{\tensor{F}} = \tensor{V}^* \tensor{F}.
\end{equation}

A popular choice of a projection subspace is a block (rational) Krylov subspace
given by
\begin{equation}
\mathcal{K}_n(\sigma) = \mbox{colspan} \left\{ 
\left( \tensor{A} + \sigma_1 \tensor{I} \right)^{-1} \tensor{F}, \ldots, 
\left( \tensor{A} + \sigma_n \tensor{I} \right)^{-1} \tensor{F} \right\},
\end{equation}
where the shifts $\sigma_j$ are distinct or repeated, finite or infinite.
Here, we use the simplest choice $\sigma_1 = \sigma_2 = \ldots = \sigma_n = 0$,
which yields a subspace
\begin{equation}
\mathcal{K}_n(0) = \mbox{colspan} \left\{ 
\tensor{A}^{-1} \tensor{F}, \tensor{A}^{-2} \tensor{F}, 
\ldots, \tensor{A}^{-n} \tensor{F} \right\},
\label{eqn:krylov0}
\end{equation} 
that can be obtained by applying a block Lanczos iteration to $(\tensor{A}^{-1}, \tensor{F})$.
Note that the computation of the basis for $\mathcal{K}_n(0)$ requires multiple 
linear solves with the matrix $\tensor{A}$. This is where the bulk of the 
computational cost of the first stage originates. However, it is alleviated by 
several factors. First, the computation is only done on small subdomains. 
Second, the computations for different subdomains are independent of each other, 
thus they can be performed in parallel. Third, the computation only must be 
done once before the time stepping. Also, a precomputed Cholesky factorization 
can be reused for the repeated linear solves. Note that unlike the ROM wave 
propagation scheme of \cite[]{pereyra2008fast}, our computation of $\tensor{V}$ 
is also independent of the number or position of sources and receivers. 

Projection subspace (\ref{eqn:krylov0}) is easy to implement and it provides
good accuracy solutions, as shown in the numerical experiments below. However,
it may not be optimal in terms of the number of degrees of freedom per
wavelength and the possible improvement of the CFL conditions. Other model
reduction techniques such as time- and/or frequency-limited balanced truncation
\cite[]{gugercin2004survey} are more appropriate for these purposes. 
Integration of these approached into our framework remains a topic of 
future research.

Although the size of the reduced order spatial operator matrix 
$\widetilde{\tensor{A}}$ is much smaller than $N$, the size of 
$\tensor{A}$, in general it is a dense matrix. In contrast, being a 
discretization of a differential operator $\tensor{A}$ is typically very 
sparse. Since the time stepping involves matrix-vector multiplications with 
$\widetilde{\tensor{A}}$, the number of non-zero entries is more 
important for the computational cost than the size of the matrix. We show in 
the next section how $\widetilde{\tensor{A}}$ can be sparsified without 
affecting $\widetilde{\tensor{M}} (\omega)$. This construction also plays 
an important role in coupling the adjacent subdomains for the proper 
exchange of information at the time-stepping stage.

\section{Stage 2: time stepping}

Once the ROMs are computed for all the subdomains, the time
stepping can be performed with the reduced order spatial operators 
$\widetilde{\tensor{A}}$. To formulate the coupling conditions for the 
adjacent subdomains and also to sparsify $\widetilde{\tensor{A}}$ we 
transform them to a special block tridiagonal form. 

Transformation to the block tridiagonal form can be done by applying a block
version of the Lanczos iteration to the pair 
$(\widetilde{\tensor{A}}$, $\widetilde{\tensor{F}})$ to obtain a unitary 
$\tensor{Q} \in \mathbb{R}^{6mn \times 6mn}$ such that 
\begin{equation}
\tensor{T} = \tensor{Q}^* \widetilde{\tensor{A}} \tensor{Q}, \quad
\tensor{R} = \tensor{Q}^* \widetilde{\tensor{F}} = [\tensor{B}_1, 0, 0, \ldots, 0 ]^*,
\end{equation}
where $\tensor{T}$ is a Hermitian block tridiagonal matix with Hermitian 
blocks $\tensor{A}_j \in \mathbb{R}^{6m \times 6m}$ on the main diagonal and 
$\tensor{B}_j \in \mathbb{R}^{6m \times 6m}$ on super/sub-diagonals.
Unitarity of $\tensor{Q}$ guarantees that the transformed transfer function
\begin{equation}
\widetilde{\tensor{M}}(\omega) = \tensor{R}^* 
\left( \tensor{T} + \omega^2 \tensor{I} \right)^{-1} \tensor{R},
\label{eqn:NtDtri}
\end{equation}
is exactly the same as (\ref{eqn:NtDROM}).

An alternative expression is available for (\ref{eqn:NtDtri}) that makes
apparent the connection to finite-difference schemes. If we apply the unitary
transformation $\tensor{V}\tensor{Q}$ to (\ref{eqn:Helmholtz}) then taking 
into account the tridiagonal structure of $\tensor{T}$, we can write
\begin{equation}
\begin{split}
\tensor{A}_1 \tensor{W}_1 + 
\tensor{B}_2 \tensor{W}_2 + \omega^2 \tensor{W}_1 & = \tensor{B}_1, \\
\tensor{B}_j \tensor{W}_{j-1} + 
\tensor{A}_j \tensor{W}_j + 
\tensor{B}_{j+1} \tensor{W}_{j+1} + \omega^2 \tensor{W}_j & = \tensor{0},
\end{split}
\label{eqn:tridiagW}
\end{equation}
for matrices $\tensor{W}_j \in \mathbb{R}^{6m \times 6m}$, $j=1,2,\ldots, n+1$, 
with $\tensor{W}_{n+1} = 0$. Then using the structure of 
$\tensor{R}$, the expression for the transfer function (\ref{eqn:NtDtri}) is simply
\begin{equation}
\widetilde{\tensor{M}}(\omega)  = \tensor{B}_1 \tensor{W}_1.
\label{eqn:NtDW}
\end{equation}

A second change of coordinates can simplify (\ref{eqn:NtDW}) even further. We can  
transform (\ref{eqn:tridiagW}) to
\begin{equation}
\begin{split}
\bhG_1 \left( \bG_1 (\tensor{U}_2 - \tensor{U}_1) \right) + 
\omega^2 \tensor{U}_1 & = \bhG_1, \\
\bhG_j \left( \bG_j (\tensor{U}_{j+1} - \tensor{U}_j) -
\bG_{j-1} (\tensor{U}_{j} - \tensor{U}_{j-1}) \right) +
\omega^2 \tensor{U}_j & = \tensor{0},
\end{split}
\label{eqn:GammaU}
\end{equation}
where $\tensor{U}_j \in \mathbb{R}^{6m \times 6m}$, $j=1,2,\ldots, n+1$, 
with $\tensor{U}_{n+1} = 0$. The corresponding transformation is done 
recursively
\begin{equation}
\begin{split}
\tensor{G}_{j+1} & = \left[ \tensor{G}_j^{*} \bG_j \right]^{-1} \tensor{B}_{j+1}, \\
\bhG_{j+1} & = \tensor{G}_{j+1} \tensor{G}^{*}_{j+1}, \\
\bG_{j+1} & = - \tensor{G}^{-*} _{j+1} \tensor{A}_{j+1} \tensor{G}^{-1}_{j+1} - \bG_j, \\
\tensor{U}_{j+1} & = \tensor{G}_{j+1} \tensor{W}_{j+1},
\end{split}
\label{eqn:GammaTrans}
\end{equation}
starting with $\tensor{G}_1 = \tensor{B}_1$, $\bG_0 = 0$. The transfer function 
is trivial
\begin{equation}
\widetilde{\tensor{M}}(\omega)  = \tensor{U}_1.
\end{equation}

In 1D, all the quantities in (\ref{eqn:GammaU})--(\ref{eqn:GammaTrans}) would
be scalars, so the following expression for the transfer function is known to be valid
\begin{equation}
\widetilde{M}(\omega) = \frac{1}{\widehat{\Gamma}_1^{-1} \omega^2 + \dfrac{1}{\Gamma_1^{-1} +
    \dfrac{1}{\ddots \; + \dfrac{1}{\widehat{\Gamma}^{-1}_n \omega^2 + \Gamma_n}}}},
\end{equation}
which is known as a Stieltjes continued fraction (S-fraction). Thus, our method
expresses the 3D NtD map as a matrix generalization of the S-fraction.
Note that for a uniform medium in 1D the scalars $\Gamma_j^{-1}$, 
$\widehat{\Gamma}_j^{-1}$ are the grid steps of a finite-difference scheme 
(\ref{eqn:GammaU}) on an optimal (spectrally matched) grid.

Relations (\ref{eqn:GammaU}) provide an easy way to obtain the coupling
conditions for the two adjacent subdomains $\Omega^\alpha$ and $\Omega^\beta$.
Let us denote by $U_j^\alpha, U_j^\beta \in \mathbb{R}^{6m}$, 
$j=1,\ldots,n$ the solution vectors on all $n$ ``layers'' of the ROM. 
The vectors $U_j$ are related to the solutions $u$ of the original equation 
(\ref{eqn:wave}) by a combined transformation 
\begin{equation}
U_j = \tensor{G}_j \left[\tensor{Q}^* \tensor{V}^* u\right]_j.
\label{eqn:Uj}
\end{equation}

To obtain time stepping for the boundary solutions $U_1^\alpha$
we match the solutions and normal fluxes on $\mathcal{B}^{\alpha\beta}$
similarly to finite-volume type methods. These matching conditions 
applied to (\ref{eqn:GammaU}) imply
\begin{equation}
\left\{
\begin{aligned}
&\frac{d^2}{d t^2} \left(  [ (\bhG_1^\alpha)^{-1} U_1^\alpha ]_\beta +
[ (\bhG_1^\beta)^{-1} U_1^\beta ]_\alpha \right) = \\
&\qquad [ \bG_1^\alpha (U_2^\alpha - U_1^\alpha) ]_\beta + 
[ \bG_1^\beta (U_2^\beta - U_1^\beta) ]_\alpha \\
&[U_1^\alpha]_\beta = [U_1^\beta]_\alpha
\end{aligned}
\right.,
\label{eqn:mass}
\end{equation}
where $[X^\alpha]_\beta$ denotes the restriction of $X^\alpha$ on
$\mathcal{B}^{\alpha\beta}$.

Note that unless $\bhG_1^\alpha, \bhG_1^\beta \in \mathbb{R}^{6m \times 6m}$ are block 
diagonal with $m \times m$ blocks, equations (\ref{eqn:mass}) define a
time-stepping scheme for the boundary solutions $U_1$ with a global mass matrix. 
This can be avoided by ensuring that the boundary functions on 
$\mathcal{B}^{\alpha \beta}$ do not overlap for all $\beta \in \mathcal{N}(\alpha)$ 
and also by adding $\tensor{F}^\alpha$ to the projection subspace. 
Then for the shared boundary solution 
$U^{\alpha\beta}_1 = [U_1^\alpha]_\beta = [U_1^\beta]_\alpha$
relations (\ref{eqn:mass}) decouple into a scheme
\begin{equation}
\begin{split}
\frac{d^2 U^{\alpha\beta}_1}{dt^2} = &
\left( [ \bhG_1^\alpha ]_\beta^{-1} + [ \bhG_1^\beta ]_\alpha^{-1} \right)^{-1} \times \\
& \times \left(
[ \bG_1^\alpha (U_2^\alpha - U_1^\alpha) ]_\beta + 
[ \bG_1^\beta (U_2^\beta - U_1^\beta) ]_\alpha
\right),
\end{split}
\label{eqn:bdrytstep}
\end{equation}
which only requires communication between the adjacent subdomains.

The time stepping for the interior solutions $U_j$, $j=2,\ldots,n$
is always fully local
\begin{equation}
\frac{d^2 U_j}{dt^2} = \bhG_j \left( \bG_j (U_{j+1} - U_j) -
\bG_{j-1} (U_{j} - U_{j-1}) \right).
\label{eqn:intrtstep}
\end{equation}

Any standard time stepping scheme can be used for 
(\ref{eqn:bdrytstep})--(\ref{eqn:intrtstep}) including Virieux,
Runge-Kutta, etc. The expressions on the right hand side of 
(\ref{eqn:bdrytstep})--(\ref{eqn:intrtstep}) are always evaluated
at the current time step.

\section{Method summary}

We summarize below our method as an algorithm that is well suited for 
parallel HPC platforms.

\begin{alg}[S-fraction multiscale finite-volume method]~\\
Stage 1. In full parallel mode for each subdomain $\Omega^\alpha$ do
the following:
\begin{enumerate}
\item[(1.1)] Compute the projection subspace bases $\tensor{V}^\alpha$
 and the reduced order models $(\widetilde{\tensor{A}}^\alpha, 
\widetilde{\tensor{F}}^\alpha)$. 
\item[(1.2)] Apply the block Lanczos algorithm to transform 
$(\widetilde{\tensor{A}}^\alpha, \widetilde{\tensor{F}}^\alpha)$ 
to a block tridiagonal form
$(\tensor{T}^\alpha, \tensor{R}^\alpha)$.
\item[(1.3)] Obtain the S-fraction coefficients $\bG_j^\alpha$, $\bhG_j^\alpha$ 
from $(\tensor{T}^\alpha, \tensor{R}^\alpha)$ using relations (\ref{eqn:GammaTrans})
\footnote{Steps (1.2) and (1.3) can be combined using a particular form of 
block Lanczos method.}.
\item[(1.4)] Project the initial conditions $u^\alpha|_{t=0}$ and 
$\partial_t u^\alpha|_{t=0}$ on the ROM subspace (\ref{eqn:Uj})
to obtain the initial conditions for $U_j^\alpha$.
\end{enumerate}
Stage 2. Starting with initial conditions $U_j^\alpha|_{t=0}$ and
$\partial_t U_j^\alpha|_{t=0}$ for each time step do the following:
\begin{enumerate}
\item[(2.1)] Exchange $[ \bG_1^\alpha (U_2^\alpha - U_1^\alpha) ]_\beta$ 
and $[ \bG_1^\beta (U_2^\beta - U_1^\beta) ]_\alpha$ between the subdomains 
sharing $\mathcal{B}^{\alpha\beta}$.
\item[(2.2)] While waiting for the data exchange, compute in parallel 
for each $\Omega^\alpha$ the updates to the interior solutions 
$U_2^\alpha,\ldots,U_n^\alpha$ using (\ref{eqn:intrtstep}).
\item[(2.3)] Once the data exchange is complete, compute in parallel 
for each $\Omega^\alpha$ the updates to the boundary solutions 
$U_1^\alpha$ using (\ref{eqn:bdrytstep}).
\end{enumerate}
\end{alg}

Note that the order of steps (2.1) and (2.2) allows for what is known in 
computer science literature as hiding the communication latency behind the 
computations. Also, observe that the communication cost is very low. We only 
exchange vectors of size $m$ between the adjacent subdomains as if we had
a second order scheme. In practice the number $m$ of boundary basis functions 
is chosen based on the source frequency and thus the minimal wavelength of 
the resulting wavefield. It does not depend on the accuracy of the reference 
fine grid discretization. This constrasts sharply with the traditional domain 
decomposition approaches for high order finite-difference schemes, where the 
communication cost is proportional to the size of the stencil. Such a small 
communication cost is possible because the ROMs approximate the NtD map to 
high (spectral) accuracy, even though (\ref{eqn:GammaU}) resemble a 
three-point difference scheme.

\section{Numerical experiments}

We study the viability of our method on a simple numerical example below. 
We consider an acoustic wave equation
\begin{equation}
u_{tt} = c^2 \Delta u,
\end{equation}
in a 3D box $\Omega = [0,7] \times [0,7] \times [0,3]$, which is split into 
$7 \times 7 \times 3$ unit cube subdomains each containing $20 \times 20 \times 20$
reference fine grid nodes. The sizes of ROMs on each $\Omega^\alpha$
are $m=25$, $n=3$.

The sound speed $c(x,y,z)$ does not depend on $z$, its dependence on $x$ 
and $y$ is shown in Figure \ref{fig:soundspeed-cskew5xymc7seg}. All the 
quantities in the example are dimensionless. First order absorbing boundary 
conditions are enforced on $\partial \Omega$. 

The numerical experiment is designed to emphasize the fact that our 
method allows for the arbitrary placement of the subdomain boundaries 
relative to the discontinuities of the coefficients of the wave equation. 
Many subdomains contain one or more discontinuity interfaces of 
$c(x,y,z)$ including corners. A thin slow fracture of contrast 
$\max(c)/\min(c) = 10/3$ passes through $\Omega$.

We simulate a single source located at $(3.5, 1.5, 1.5)$ that emits a
Gaussian pulse corresponding to a minimal wavelength of $\lambda=0.78$ for 
$c=1$. The solution traces $d(x,t) = u(x, 0, 1.5, t)$ are measured for 
$t \in [0, 12.5]$. For easier visualization, we normalize the traces by 
$\int_0^7 d(x,t) dx$ for each $t \in [0.73,8.2]$ with the results given in 
Figure \ref{fig:seismonorm_fine_cskew5xymc7seg,seismonorm_rom_cskew5xymc7seg}.

We observe in Figure \ref{fig:seismonorm_fine_cskew5xymc7seg,seismonorm_rom_cskew5xymc7seg}
a good agreement between our method's solution and the solution obtained on
a reference fine grid. The relative $L_2$ norm error between the two is $2.7\%$.
This is achieved with six reduced order degrees of freedom per wavelength per 
dimension compared to $16$ points per wavelength for the fine grid scheme.

\plot{soundspeed-cskew5xymc7seg}{width=0.9\columnwidth}
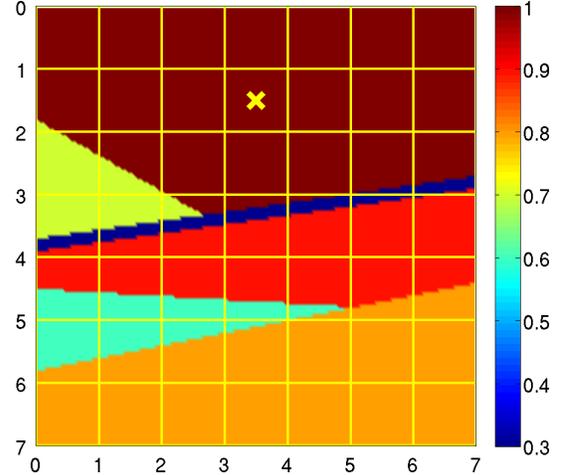
{Sound speed profile at $z=1.5$. Subdomain boundaries are yellow lines, 
source location is at $\times$.}

\multiplot{2}{seismonorm_fine_cskew5xymc7seg,seismonorm_rom_cskew5xymc7seg}{width=\columnwidth}
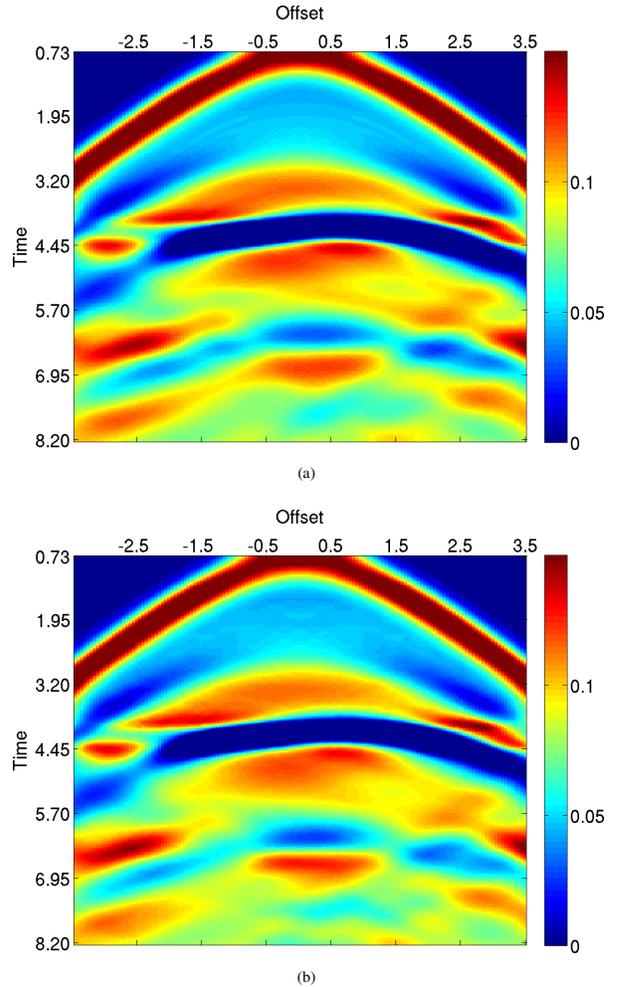
{Solution traces at $y=0$, $z=1.5$ normalized for each time for easier visualization:
(a) reference fine grid; (b) S-fraction multiscale finite-volume method. Relative error 
for for $t \in [0, 12.5]$ is $2.7\%$.}

\section{Conclusions and future work}

We performed a first study of a general framework for numerical wave propagation
in the time domain using the ROMs in a multiscale setting. In the early numerical 
experiments the method demonstrated good accuracy while substantially reducing the 
number of degrees of freedom per wavelength compared to the traditional FDTD schemes.
Further improvements in the model order reduction should allow us to approach the 
Nyquist limit both in space (fewer points per wavelength) and in time 
(relaxed CFL conditions).


\twocolumn

\bibliographystyle{seg}  
\bibliography{segabs_sfmsfv}

\end{document}